\documentclass[12pt]{article}
\usepackage{amsmath,amsthm,amssymb}

\newtheorem{theorem}{Theorem}

\newtheorem{lemma}{Lemma}
\newtheorem{definition}{Definition}
\newtheorem{corollary}{Corollary}

\newcommand{\gd}{\delta}
\newcommand{\gl}{\lambda}
\newcommand{\gs}{\sigma}
\newcommand{\hess}{\nabla^2}
\newcommand{\gra}{\nabla}
\newcommand{\de}{\partial}
\newcommand{\bpf}{\begin{proof}}
\newcommand{\epf}{\end{proof}}
\newcommand{\beq}{\begin{equation}}
\newcommand{\eeq}{\end{equation}}

\begin{document}
\title{ Local Estimates for Some Fully Nonlinear Elliptic Equations}
 \vskip 1em
 \author{Szu-yu Sophie Chen}
 \date{September 15, 2005}
\maketitle

 \begin{abstract}
 We present a method to derive local estimates for some classes of fully
 nonlinear elliptic equations. The advantage of our method is that
 we derive Hessian estimates directly from $C^0$ estimates.
 Also, the method is flexible and can be applied to a large class
 of equations.
 \end{abstract}
\vskip 1em

 Let $(M,g)$ be a smooth Riemannian manifold of dimension
$n \geq 2$. We are interested in a priori estimates for solutions
of some classes of fully nonlinear elliptic equations on $(M, g).$
These kinds of equations arise naturally from geometry and other
fields of analysis and share structures similar to those of the
Monge-Ampere equations.

Regularity problems are studied by people in different fields
separately. One would like to ask if it is possible to give a
unified proof and to generalize further to a large class of
equations. The answer is affirmative provided the equations
satisfy some algebraic structures which can induce the
cancellation phenomenon. We will see how this phenomenon helps us
to get the Hessian bound directly. One of the interesting cases is
the Schouten tensor equation arising from conformal geometry:
 $$\gs_k ^{\frac{1}{k}} (g^{-1}(\nabla ^2 u + du \otimes
du -\frac{1}{2}|\nabla u|^2 g + A_g) ) = f(x)\, e^{-2u}$$
 where $\gs_k$ is the kth elementary symmetric function. Local $C^2$
 estimates are proved for this equation by Chang, Gursky, and Yang \cite{CGY02}
 ($k=2,n=4$) and by P.Guan and G.Wang \cite{GW03} for all $k \leq n$.
  The same results with specific dependence on the radius of the
  domain are established by Gursky and Viaclovsky \cite{GV05}.
 P.Guan and G.Wang \cite{GW04} also prove the local estimates for
 quotients of the elementary symmetric
 functions. Other related works in this direction include
 \cite{LL03}, \cite{GV03} and \cite{STW05}.

 Another interesting case is the following equation in optics
 geometry:
$$\det (g_c^{-1} (\hess v - \frac{|\gra v|^2}{2v} g_c + \frac{v}{2} g_c)) =
 \left(\frac{|\gra v|^2 + v^2}{2v}\right )^n \nu(x) \phi(S(x, v, \gra v)). $$
  Interior $C^2$ estimates are proved by X.Wang \cite{Wx96} for $n=2,$
  while local $C^2$ estimates for all $n \geq 2$ are by P.Guan and X.Wang \cite{GWx98}.

 It turns out that in getting local $C^2$ estimates for nonlinear equations as described
 above, the coefficient in front of the gradient square
 term plays an important role. For arbitrary coefficients, in general it is not true
 that we have local estimates . See \cite{STW05} for a counterexample.

 In the degenerate case, when the gradient square term disappears,
 one can have maximal principles for second derivatives, which means
 the Hessian bound over $\overline{\Omega}$ is less than or equal to
 that on $\de \Omega.$ Examples are the general Monge-Ampere
 equations. In particular, the Gauss curvature equation in
 a domain $\Omega$ in $\mathbb{R}^n,$
 $$ \det^{\qquad \frac{1}{n}}(\hess u)= \kappa^{\frac{1}{n}} (x)
  ( 1 +|\gra u|^2)^{\frac{n+2}{2n}},$$
  is of this type. Another relevant equation is the Gauss curvature
  equation for a radial graph in a domain $\Omega$ in $S^n$:
$$\det^{\qquad \frac{1}{n}}(g_c ^{-1} (\hess v + v g_c))
  = \kappa^{\frac{1}{n}} (x) ( \frac{v^2 + |\gra v|^2}{v^2})^{\frac{n+2}{2n}}.$$
 Maximal principles for the first equation is studied in
 Caffarelli, Nirenberg, and Spruck \cite{CNS-1}, and
 for the latter by B.Guan and J. Spruck \cite{GbS93}. See also \cite{gbs02},
 \cite{TrWx02} and \cite{GGb02} for related works.

 In this paper, we consider more general operators, which
  in particular include the equations discussed above.
 We will derive local $C^2$ estimates directly
 from $C^0$ bounds and also prove the maximal principles for second
 derivatives. For the reader who is more interested in the
 aforementioned equations, he or she can jump directly
 to Section~\ref{s:ex} where a brief explanation and statements
 of results about them are given.

  Now we turn to the equations we are going to discuss.
 Let \beq \label{e:w}
 W= \hess u + a(x) du \otimes du + b(x) |\gra u|^2 g +
 B(x) \eeq
  be a $(0,2)$ tensor on a Riemannian manifold $(M^n, g)$.
  The derivatives are covariant derivatives with respect to the
  metric $g$. Consider the equation
  \beq \label{e:main}
 F(g^{-1}W)= f(x,u) h(x,\gra u)
 \eeq where $F$ satisfies some fundamental structure conditions
 listed later and  $g^{-1}$ is the induced inverse tensor of
 metric tensor $g.$ Equation (\ref{e:main}) means that
 we apply $F$ to the eigenvalues of matrix (or $(1,1)$ tensor)
 $g^{-1} W.$ When the manifold is flat (e.g., the Euclidean case), we have
 $g_{ij} = \gd_{ij}$ where $\gd_{ij}$ is the Kronecker delta.
 In this case, we drop $g^{-1}$ and simply write
 $F(W)= f(x,u) h(x,\gra u).$

  We now describe the fundamental structure conditions for $F.$
 Let $\Gamma$ be an open convex cone
 with vertex at the origin satisfying $\{\gl: \gl_i > 0 , \forall i \} \subset \Gamma
 \subset \{\gl: \sum_i \gl_i > 0\}.$ Suppose that $F(\gl)$ is a homogeneous
 symmetric function of degree one in $\Gamma$ normalized with $F(e) = F((1,\cdots,1)) = 1.$
 Moreover, $F$ satisfies the following in $\Gamma:$

 (S0) $F$ is positive.\par
 (S1) $F$ is concave. (i.e., $\frac{\de^2 F}{\de \gl_i \de
 \gl_j}$ is negative semi-definite.)\par
 (S2) $F$ is monotone. (i.e., $\frac{\de F}{\de \gl_i}$ is
 positive.)\par

  In some cases, we need an additional condition:\par
 (A) $\sum_i \frac{\de F}{\de \gl_i} \geq \mu_0 \left( \frac{\sum_i
 \gl_i}{F}\right)^{\mu_1}$, for some $\mu_0, \mu_1 > 0.$

 \vskip 1em

 An easy example is $F = \gl_1 + \cdots + \gl_n$ with $\Gamma =
 \{\gl : \gl_1 + \cdots + \gl_n > 0 \}.$ Then $F(g^{-1} W) = tr_g\, W =
 \Delta u + (a(x)+ nb(x)) |\gra u|^2 + tr_g\, B(x) $ is just
 the Laplace-Beltrami operator plus some lower order terms, where $tr_g$
 is the trace with respect to $g$. More interesting examples are discussed in
 Section~\ref{s:sigma}. Condition (S1) is necessary in most
 elliptic theories. Condition (S2) is the actual ellipticity. It
 is an elementary fact that if $F$ is a symmetric function of
 eigenvalues, then $\frac{\de F}{\de \gl_i}> 0 $ for all $i$ if and only if
 $F^{ij} :\equiv \frac{\de F}{\de W_{ij}}$ is positive definite.
 Condition (A) is used previously in \cite{GWx98}.
 There is one more key point. In general, we \emph{do not} have uniform
 ellipticity for fully nonlinear elliptic equations. This is because
 $F^{ij}$ involves $\hess u$ whose a priori estimates need to be
 derived.

 A natural question is whether we can consider
 the tensor in forms other than (\ref{e:w}). It turns out that
 for some equations coming from geometry, they can be formulated
 in the form of (\ref{e:w}) after a wise choice of the function $u(x).$
 We will see in Section~\ref{s:ex} that it is certainly the case
 for geometric optics equations and Gauss curvature equations on spheres.

 Before stating the theorems, we introduce the following notations. Let\\
 $f(x,z) : M^n \times \mathbb R \rightarrow \mathbb R$ and
 $h(x,p): M^n \times \mathbb{R}^n \rightarrow \mathbb R$
 be two given positive functions.\\ Let $u= u(x) : M^n \rightarrow \mathbb R$
 be a solution to (\ref{e:main}). We define
 \begin{align}
 c_{inf} =& \inf_{x \in M} f(x,u) , \notag\\
 c_{sup} =& \sup_{x \in M} ( f+ | \gra_x f(x, u)|+ |f_z (x, u)|+
    | \hess_x f(x, u)| + |\gra_x f_z (x, u)| +
 |f_{z z} (x, u)|) ,\notag\\
 e_{sup} =& \sup_{x \in M} ( h(x, \gra u) + | \gra_x h (x, \gra u)|+ |\gra_p h (x,\gra
 u)|+ | \hess_p h(x, \gra u)| \notag \\
         &   + |\gra_x \gra_p h (x, \gra u)| +
          |\hess_x h (x, \gra u)|).\notag
 \end{align}
 If we restrict $x$ to a local ball $B_r$, we use the
 corresponding notations $c_{inf} (r)$, $c_{sup}(r)$ and $e_{sup}(r).$
 We also use the convention that a $(0,2)$ tensor $T_{ij} \geq g_{ij}$
 means that $T(v,v) \geq g(v,v)$ for all vectors $v$ as a bilinear
 form.

 In the following theorem, cases (a) and (b) show how relations of
 $a(x), b(x)$ and $h(x, \gra u)$ give us Hessian estimates directly.
 However, if we have gradient bounds already,
 then $h$ is bounded. Case (c) shows more general results.

 \begin{theorem}\label{t:local}(Local estimates)
  Let $F$ satisfy the structure conditions (S0)-(S2) in a corresponding cone
  $\Gamma$ and $u(x)$ be a $C^4$ solution to (\ref{e:main}) in a local geodesic ball
  $B_r.$
  Suppose that $b(x) < - \gd_1$ and $a(x) + n b(x) < -\gd_2.$\par

  \underline{case(a)}: $h= h_0$ is a positive constant.
  Then $$ \sup_{x\in B_{\frac{r}{2}}} (|\hess u| + |\gra u|^2)
   \leq C_1,$$ where $C_1 = C_1( n, r, \|a\|_{C^2}, \|b\|_{C^2},
   \|B\|_{C^2}, \|g\|_{C^3}, h_0, \gd_1, \gd_2, c_{sup}(r))$
   but is independent of $c_{\inf}(r).$\par
  \underline{case(b)}: Suppose $h = h(\gra u)$ and $f=
  f(x)$. Let $\Lambda(p)$ be a positive function such that $h_{p_i p_j}\geq \Lambda(p) g_{ij}.$
  If there exists some number $M > 0$ such that
  $$ h \leq M \Lambda(p) (1 + |p|)^2 \qquad and \qquad  |\gra_p h| \leq M \Lambda(p) (1+ |p|),$$
  then
  $$ \sup_{x\in B_{\frac{r}{2}}} (|\hess u| + |\gra u|^2) \leq C_2,$$ where
  $C_2 = C_2(n, r, \|a\|_{C^2}, \|b\|_{C^2}, \|B\|_{C^2}, \|g\|_{C^3},
  \gd_1, \gd_2, M,\sup_p \Lambda(p), c_{sup}(r), c_{inf}(r)).$\par
  \underline{case(c)}: Suppose that $F$ satisfies the additional
  condition (A) and that $\Gamma^+_2 \subset \Gamma$. ( See Section~\ref{s:sigma} for the
   definition of $\Gamma^+_2.$)  Then
    $$ \sup_{x\in B_{\frac{r}{2}}} |\hess u| \leq C_3,$$ where
    $C_3 = \: C_3 (\mu_0, \mu_1, n, r, \|a\|_{C^2}, \|b\|_{C^2}, \|B\|_{C^2}, \|g\|_{C^3},
    \gd_1, \gd_2, e_{sup} (r), c_{sup} (r)$,\\ $\sup_{B_r} |\gra u| ).$
 \end{theorem}

 An example of case (a) is the Schouten tensor equation
 arising from conformal geometry; an example of case (b) and (c)
 is the geometric optics equation.

 For the degenerate case $b= 0$, we do not in general have
 local estimates. However, if the manifold has enough symmetry, say
 of constant sectional curvature $K$, we may consider a special type
 of equation
 \beq \label{e:csc}
 F(g^{-1}(\hess u + a du \otimes du + K g )) = f(x,u) h(\gra u)
 \eeq where $a$ is a constant. Note that when $a = 0$ and $K= 0$
 (e.g., Euclidean space), this is the Monge-Ampere type equation.

 \begin{theorem}\label{t:max}(Maximum Principle)
  Let $F$ satisfy the structure conditions (S0)-(S2) in a corresponding
   cone $\Gamma$. Suppose that $(M,g)$ is of nonnegative constant sectional
   curvature $K$ and that $h_{p_i p_j} \geq \epsilon \gd_{ij}$ for
   some positive $\epsilon.$ Let $u(x)$ be a $C^4$ solution to (\ref{e:csc})
   in a bounded domain $\Omega \subset M.$
  Then
   $$\sup_{x \in \bar{\Omega}} |\hess u| \leq C_4$$
  where $C_4 = C_4 ( n, a, K, \epsilon, e_{sup}, c_{sup}, c_{inf}, \|u\|_{C^1(\overline{\Omega})},
  \sup_{\de \Omega} |\hess u| ).$
 \end{theorem}

 Examples of Theorem~\ref{t:max} are the Gauss curvature equations on a domain in $R^n$
 and in $S^n.$

 This paper is organized as follows. We start with some background in
 Section~\ref{s:sigma}. In Section~\ref{s:ex}, we discuss
 applications and give the statements of results. The proofs of
 Theorems~\ref{t:local} and \ref{t:max} are in Sections~\ref{s:local} and
 \ref{s:max}, respectively.
 \vskip 1em

 \textbf{Acknowledgments:} The author would like to thank Alice
 Chang for drawing the author's attention to the work
 about Calabi-Yau problem, where $C^2$ bounds are
 derived directly. Otherwise, this paper could not
 be possible. The author also appreciates Matt Gursky's
 valuable and enlightening math discussions.

\section{Background}\label{s:sigma}

 First, we give some basic facts about homogeneous symmetric functions.
\begin{lemma}\label{l:sym}
 Let $\Gamma$ be an open convex cone with vertex at the origin
 satisfying $\{\gl: \gl_i > 0, \forall i\} \subset \Gamma$ and
 $e = (1, \cdots, 1)$ be the identity.
 Suppose that $F$ is a homogeneous symmetric function of degree one
 normalized with $F(e)= 1,$ and that $F$ is concave in $\Gamma.$
 Then the following are true:
 \begin{description}
 \item{(a)} $\sum_i \gl_i \frac{\de F(\gl)}{\de \gl_i} = F(\gl), \quad$ for $\gl \in
 \Gamma.$
 \item{(b)} $\sum_i \frac{\de F(\gl)}{\de \gl_i} \geq F(e) = 1, \quad$
 for $\gl \in \Gamma.$
 \end{description}
\end{lemma}
 \bpf (a) By homogeneity, $F(\theta \gl) = \theta F(\gl).$
  Let $\theta$ be some positive number. Since $F$ is concave
  in $\Gamma$, then
 $$ (\theta - 1) F(\gl) = F(\theta \gl) - F(\gl) \leq
 \sum_i (\theta \gl_i - \gl_i) \frac{\de F(\gl)}{\de \gl_i}.$$
 Choose some $\theta < 1$ and some $\theta > 1$ and cancel out the
 factor $(\theta -1),$ which gives the result.\par
  (b) $\Gamma$ contains the identity $e$ and since $F$ is concave in
    $\Gamma$, we have
    $$ F(e) - F(\gl) \leq \sum_i (1 - \gl_i) \frac{\de F(\gl)}{\de \gl_i} =
    \sum_i \frac{\de F(\gl)}{\de \gl_i} - F(\gl)$$
  where the equality holds by (a). Cancelling out $F(\gl)$, we prove (b).
 \epf

 Now we focus on elementary symmetric functions because most interesting
 cases are related to them.

\begin{definition} Let $W$ be a matrix with eigenvalues $\gl_1,\cdots,\gl_n.$
 Then $\gs_k (\gl(W))= \sum_{i_1<\cdots< i_k}
\gl_{i_1}\gl_{i_2}\cdots\gl_{i_k}$ for $k \leq n$ is called the
kth elementary symmetric function of the eigenvalues of $W$. We
denote $\gs_0 = 1.$ For examples, $\gs_1 = \gl_1 + \cdots + \gl_n
= tr \, W$ and $\gs_n = \gl_1 \cdots \gl_n = \det W.$
\end{definition}
 The elementary symmetric functions are examples of hyperbolic
 polynomials introduced by Garding \cite{Gar59} ,which have nice
 properties in the associated cones.
\begin{definition}
 The set $\Gamma^+_k = \{$ the connected component of $\gs_k(\gl)>
0$ which contains the identity $\}$ is called the positive k-cone.
 Equivalently, it is shown in \cite{Gar59} that for $k>0$,
 $\Gamma^+_k = \{ \gl \,: \gs _i (\gl)> 0 , 1 \leq i \leq k \}$
 is an open convex cone with vertex at the origin, e.g.,
 $\Gamma^+_1 = \{\gl : \gl_1 + \cdots \gl_n > 0\}$ and $\Gamma^+_n =
 \{\gl : \gl_i > 0 , 1 \leq i \leq n \}.$ We also have the nested
 relation
  $$ \Gamma^+_1 \supset \Gamma^+_2 \supset \cdots \supset \Gamma^+_n.$$
 We say that $W \in \Gamma^+_k$ if the eigenvalues $\gl(W) \in \Gamma^+_k.$
\end{definition}

 We list some basic properties of elementary symmetric
 functions.
\begin{lemma} \label{l:sigma}(see \cite{Gar59}, \cite{Mit70} and \cite{Tr90} for the proof)
  Let $G =(\frac{\gs_k}{\gs_l})^{\frac{1}{k-l}}, 0 \leq l < k \leq
  n.$\par

  (a) $G$ is positive and concave in $\Gamma^+_k.$\par
  (b) $G$ is monotone in $\Gamma^+_k,$ i.e., the matrix
  $G^{ij} = \frac {\de G}{\de W_{ij}}$ is positive definite.\par
  (c) For $1 \leq m < k \leq n,$ we have the Newton-MacLaurin inequality
    $$ k(n-m+1) \gs_{m-1} \gs_k \leq m(n-k+1) \gs_m \gs_{k-1}.$$
\end{lemma}

  Therefore, $F = \binom{n}{k}^{-\frac{1}{k-l}} \binom{n}{l}^{\frac{1}{k-l}}\, G$
  satisfies the structure conditions (S0)-(S2) in $\Gamma^+_k.$
  We further show that if $l = 0$ and $k \geq 2,$ then
  $F = \binom{n}{k} ^{- \frac{1}{k}} \gs_k^{\frac{1}{k}}$
  satisfies (A) with $\mu_0 = n^{-\frac{1}{k-1}}$ and $\mu_1 = \frac{1}{k-1}$.

  By Lemma~\ref{l:sigma} (c), for $1 \leq m \leq k-1,$ we have the
  recursive formula
   $$ \gs_m \geq \frac{k(n-m+1)}{m(n-k+1)} \left(\frac{\gs_{k}}{\gs_{k-1}}\right) \gs_{m-1}.$$
  Then
  $$\gs_{k-1}
    \geq \frac{k^{k-2}(n-k+2)\cdots(n-1)}{(n-k+1)^{k-2}(k-1)!}
    \left(\frac{\gs_k}{\gs_{k-1}}\right)^{k-2} \gs_1 =
    \frac{\binom{n}{k-1} ^{k-1} }{n \binom{n}{k} ^ {k-2}}
    \left(\frac{\gs_k}{\gs_{k-1}}\right)^{k-2} \gs_1,$$
  which implies
  $$ \sum_i \frac{\de F}{\de \gl_i} = \binom{n}{k} ^{- \frac{1}{k}}
    \frac{n-k+1}{k} \gs_{k-1} \gs_k ^{-\frac{k-1}{k}}
     \geq n^{-\frac{1}{k-1}}
    \left(\frac{\gs_1}{F}\right)^{\frac{1}{k-1}}.$$

  Another useful function, which is also a variant of the elementary
  symmetric functions, is
  $$ \gs_k^{\frac{1}{k}}( t \gl + s \gs_1(\gl) e).$$
   Suppose $t,s \geq 0$ with $t+s \geq 1.$
   Let $\Gamma = \{\gl: t \gl + s \gs_1(\gl) e \in \Gamma^+_k\}.$
   Then $\Gamma$ is an open convex cone with vertex at the origin.
   Let $F = \frac{1}{t+ ns} \binom{n}{k}^{-\frac{1}{k}}
   \gs_k^{\frac{1}{k}}( t \gl + s \gs_1(\gl) e).$ It is easy to
   see that $F$ is a homogeneous symmetric function of degree one.
   Moreover, it is shown in \cite{LL03} that $F$ is concave in $\Gamma.$

   Since we consider equations on manifolds, all derivatives
   are the covariant derivatives with respect to the metric $g.$
   Let $u$ be a function on a manifold. Recall that $u_{ij} = u_{ji}.$
   However, when we consider higher order derivatives, we should get some
   curvature terms if we change the order of differentiations.
   We denote the Riemannian, Ricci, and scalar curvature by $R_{ijkl}, R_{ij}$
   and $R$, respectively. The following formulae are very useful. We remind
   the readers that we assume $g_{ij}(0) = \gd_{ij}$ without loss of
   generality:
 $$\begin{array} {l}
 u_{kij} = u_{ijk} + R_{mikj} u_m\\
 u_{ijkl} = u_{ijlk} + R_{mjkl} u_{mi}+ R_{mikl} u_{mj} \\
 u_{kkij} = u_{ijkk} + 2 R_{mikj} u_{mk} - R_{mj} u_{mi} - R_{mi} u_{mj} - R_{mi,j} u_m
      + R_{mikj,k} u_m
 \end{array}$$
 Hence,
 $$\begin{array} {l}
 u_{kij} = u_{ijk} + O(|\gra u|)\\
 u_{kkij} = u_{ijkk} + O(|\hess u| + |\gra u|).
 \end{array}$$

 \section{Applications} \label{s:ex}
  In this section, we will list examples where Theorem~\ref{t:local} and
  Theorem~\ref{t:max} can be applied.
  \vskip 1em
 1.  \emph{Schouten tensor and conformal geometry}
  \vskip 1em
 Let $(M,g)$ be a smooth compact Riemannian manifold of dimension
$n \geq 3$. The Schouten tensor of $g$ is defined as
$$ A_g =\frac{1}{n-2} ( Ric - \frac{R}{2(n-1)} g ). $$
Under the conformal change $g_u= e^{-2u} g$, the tensor $A_{g_u}$
 satisfies
 $$ A_{g_u} = \nabla ^2 u + du \otimes du - \frac{1}{2}|\nabla u|^2 g +
 A_g.$$ We consider the equation ($0\leq l < k \leq n$)
 \beq \label{e:quot}
 (\frac{\gs_k}{\gs_l})^{\frac{1}{k-l}}(g^{-1}(\nabla ^2 u + du \otimes du
 -\frac{1}{2}|\nabla u|^2 g + A_g) ) = f(x)\, e^{-2u}.
 \eeq
 Local estimates are proved by Chang,Gursky, and Yang \cite{CGY02}
 ($k=2, l=0$ or $1, n=4$), Guan-Wang \cite{GW03} ($l=0$), and Guan-Wang \cite{GW04}
 with the additional assumption $(n-k +1)(n-l+1) > 2(n+1).$
  As a corollary of Theorem~\ref{t:local} (a), we prove
 local $C^2$ estimates for all $0\leq l < k \leq n$ with specific
 dependence on the radius. The following argument is a modification of
 that in Gursky and Viaclovsky \cite{GV05} where the case $l= 0$
 is proved.

 \begin{corollary} \label{c:quot}
  Let $u(x)$ be a $C^4$ solution to (\ref{e:quot}) with $A_{g_u} \in \Gamma_k^+$
  in a geodesic ball $B_r$ in $(M,g)$, $0 \leq l < k \leq n$. Suppose
  that $f$ is positive. Then
  \begin{equation} \label{i:local1}
  \sup_{B_{\frac{r}{2}}} \,( |\gra u(x)|^2 + |\hess u(x)| ) \leq
   C ( r^{-2} + \sup_{x \in B_r} e^{ -2 u})
 \end{equation}
 where $C$ depends on $ n, k, l, \|g\|_{C^4}, \| f\|_{C^2}$ but
 does not depend on $\inf f$
  \end{corollary}
 \bpf In Section~{\ref{s:sigma}} , we showed that
 $F = \binom{n}{k}^{-\frac{1}{k-l}} \binom{n}{l}^{\frac{1}{k-l}}\,
  \left(\frac{\gs_k}{\gs_l}\right) ^{\frac{1}{k-l}}$ satisfies the structure conditions
  (S0)-(S2). Let us check the conditions
  in Theorem~\ref{t:local} (a). In this case, $a = 1, b=
  -\frac{1}{2}$ and $h_0 = 1.$ Therefore, $\gd_1 = \frac{1}{2}, \gd_2 = \frac{n-2}{2}$ and $B =
  A_{g}$ whose $C^2$ norm depends on $\|g\|_{C^4}.$ Finally,
  $$c_{sup} = \binom{n}{k}^{-\frac{1}{k-l}} \binom{n}{l}^{\frac{1}{k-l}}
  \sup_{x \in B(1), i, j} (4f + |2 f_i| + |f_{ij}|) e^{-2u} \leq c \sup_{x \in B(1)} e^{-2u} $$
  where $c$ depends on $n, k, l,\| f\|_{C^2}$ but does not
  depend on $\inf f.$ Hence, we prove the case for $r= 1.$
   As for general $r,$ without loss of generality, we may assume
   the injectivity radius $\iota$ is greater or equal to one and
   $r < 1.$ Define the mapping
    $$\begin{array}{rl}
     E(y) : B_1 \subset \mathbb{R}^n \rightarrow& B_r \subset M^n\\
       y \rightarrow& \exp (r y)= x
     \end{array}$$ where $\exp$ is the exponential map. On $B(1)$,
     define the metric $\tilde g = r^{-2} E^*g$ and the function
     $\tilde u (y) = u(E(y))- \ln r.$ Then $\tilde{u}$ satisfies
     $$(\frac{\gs_k}{\gs_l})^{\frac{1}{k-l}} (\tilde{g}^{-1}(\hess_{\tilde g} \tilde u +
   d\tilde u \otimes_{\tilde g} d\tilde u -
   \frac{1}{2}|\gra_{\tilde g} \tilde u|^2 \tilde g + A_{\tilde g} )) =
   f(E(x))\, e^{-2\tilde u}$$ on $B_1.$ By the estimates we
   obtained for $r=1$, we get
   $$\sup_{B_{\frac{1}{2}}} \,( |\gra_{\tilde g} \tilde u|^2 + |\hess_{\tilde g}
   \tilde u|) (y)
    \leq C (1 + \sup_{y \in B_1} e^{-2 \tilde u}).$$
   Now by the definitions of $E$, $\tilde{g}$ and $\tilde{u},$
  it is not hard to see $ |\gra u(x)|^2 + |\hess u(x)| = r^{-2}
  ( |\gra_{\tilde g} \tilde u|^2 + |\hess_{\tilde g} \tilde u|) (y)$
  and $e^{-2 \tilde u} = r^2 e^{-2u}.$ It remains to verify the
  conditions on the constant $C$. Since $r < 1,$ we have $\|\tilde g\|_{C^4(B_1)}
   \leq \|g\|_{C^4(B_r)}$ and $\|E^*f\|_{C^2(B_1)} \leq \|f\|_{C^2(B_r)}.$
 \epf
  In \cite{LL03} and \cite{GV03}, they consider the following
  equations
   $$ \gs_k^{\frac{1}{k}}( t \gl(A_{g_u}) + s \gs_1(\gl(A_{g_u})) g) = f(x, u)$$
   for $f(x, u) = f_0 (x) e^{-2u}$ and $f(x, u) = f_0(x) e^{2u}$, respectively, with
   $ t \gl + s \gs_1(\gl) g \in \Gamma^+_k , t, s \geq 0$ and $t+s \geq 1.$
  The local estimates are derived in \cite{LL03} and \cite{GV03}
  accordingly. We reprove these results as a corollary.

 \begin{corollary} Let $f_1(x,u) = f_0 (x) e^{-2u}$ and
 $f_2(x,u) = f_0(x) e^{2u}.$ Suppose $u_i(x)$ is a $C^4$ solution of
 the following equations in a geodesic ball $B_r:$
  $$ \gs_k^{\frac{1}{k}} (t \gl(A_{g_{u_i}}) + s \gs_1(\gl(A_{g_{u_i}})) g) = f_i(x, u_i) $$
  for $i= 1$ or $2$ , $ t \gl + s \gs_1(\gl) g \in \Gamma^+_k ,
  t+s \geq 1$ and $t+ ns \leq c_0.$ Then
 \begin{equation} \label{i:local2}
  \sup_{B_{\frac{r}{2}}} \,( |\gra u_1(x)|^2 + |\hess u_1(x)| ) \leq
  C ( 1 + \sup_{x \in B_r} e^{ -2 u_1})
 \end{equation} and
\begin{equation} \label{i:local3}
  \sup_{B_{\frac{r}{2}}} \,( |\gra u_2(x)|^2 + |\hess u_2(x)| ) \leq
  C ( 1 + \sup_{x \in B_r} e^{2 u_2})
 \end{equation}
  where $C = C(n, k, r, \|g\|_{C^4},\|f_0\|_{C^2})$ but is independent
  of $t, s$ and $\inf f_0.$
 \end{corollary}
 \bpf
 The proof is similar to that of Corollary~\ref{c:quot}. Let $F(\gl) = \frac{1}{t+ ns}
 \binom{n}{k}^{-\frac{1}{k}} \gs_k^{\frac{1}{k}}( t \gl + s \gs_1(\gl)
 g),$ so $F$ satisfies (S0)-(S2).
 \epf
\vskip 1em

  2. \emph{Optics Geometry}
  \vskip 1em

  Let $(S^2, g_c)$ be the standard 2-sphere. Suppose there is a
  point source light at the origin with the density function
  $\nu(x), x \in S^2$ and the light reflects according to the geometric optics.
   Given domains $\Omega, D \subset S^2,$ we are asked to find
  a star-shaped surface $\Sigma \subset \mathbb{R}^3$ whose projection to $S^2$
  is $\Omega$ such that the light reflected from $\Sigma$ travels in directions in $D$
  with density $\phi^{-1}(x), x \in D.$ This is related to the reflector antenna design
  problem. Mathematically, it means to find a positive solution $v$ of the
  fully nonlinear elliptic equation
 $$\det (g_c^{-1} (\hess v - \frac{|\gra v|^2}{2v} g_c + \frac{v}{2} g_c)) =
 \left(\frac{|\gra v|^2 + v^2}{2v}\right )^2 \nu(x) \phi(S(x, v, \gra v)) $$
  where $S(x, v, \gra v) = - \frac{2v \gra v + (v^2 - |\gra v|^2)N(x)}{|\gra v|^2 +
  v^2}$ and $N(x)$ is the unit vector pointing to
  $x \in S^2.$ (For background and derivation of the equation, see \cite{WN75},
  \cite{Wx96} and \cite{GWx98}.) Let us consider the general equation on $S^n:$
  $$\det (g_c^{-1} (\hess v - \frac{|\gra v|^2}{2v} g_c + \frac{v}{2} g_c)) =
 \left(\frac{|\gra v|^2 + v^2}{2v}\right )^n \nu(x) \phi(S(x, v, \gra v)). $$
 The tensor inside $\det$ is not in the form of (\ref{e:w}). However, since $v$ is positive,
 let $u= \ln v$. The equation becomes
 \beq \label{e:optic}
 \det (g_c^{-1}(\hess u + du \otimes du -\frac{1}{2}|\gra u|^2 g_c + \frac{g_c}{2})) =
 \left(\frac{|\gra u|^2 + 1}{2}\right)^n \nu(x) \phi(T(x, \gra u))
  \eeq where $T(x, u) = - \frac{2 \gra u + (1- |\gra u|^2) N(x)}{1 + |\gra
  u|^2}.$ In \cite{GWx98}, local $C^2$ estimates are proved.
   As a corollary of Theorem~\ref{t:local}(c), we prove the following.
 \begin{corollary}  Let $\;u(x)$ be a $C^4$ solution to (\ref{e:optic}) in
 a geodesic ball $\;B_r$ with\\
  $ \hess u + du \otimes du -\frac{1}{2}|\gra u|^2 g_c +
  \frac{1}{2} g_c \in \Gamma^+_n.$
  Then $$\sup_{B_{\frac{r}{2}}} \, |\hess u(x)|  \leq C$$
  where C depends on $n, r, \|\nu\|_{C^2}, \|\phi\|_{C^2}, \sup_{B_r} |\gra u|$
  but does not depend on $\inf \nu$ and $\inf \phi.$
 \end{corollary}
 \bpf In Section~\ref{s:sigma}, we showed that $F = \gs_n^{\frac{1}{n}}$
 satisfies (S0)-(S2) and (A) with $\mu_0 = n^{-\frac{1}{n-1}}$ and $\mu_1 = \frac{1}{n-1}.$
 Besides, in our case, $f = \nu^{\frac{1}{n}}(x)$
  and $h = \frac{1}{2} \phi^{\frac{1}{n}}(x, \gra u)(1 + |\gra u|^2).$
 \epf
 For a special case when $\phi$ is a positive constant, we can
 prove local $C^2$ estimates without using the gradient bound.
 This is a corollary of Theorem~\ref{t:local}(b).

 \begin{corollary} Suppose that $\phi= \phi_0$ is a positive
  constant. Let $u(x)$ be a $C^4$ solution to (\ref{e:optic}) in
 a geodesic ball $B_r$ with
  $ \hess u + du \otimes du -\frac{1}{2}|\gra u|^2 g_c +
  \frac{1}{2} g_c \in \Gamma^+_n.$
  Then $$\sup_{B_{\frac{r}{2}}} \, (|\gra u|^2+ |\hess u(x)|)  \leq C$$
  where C depends on $n, r, \phi_0, \|\nu\|_{C^2}$ and $\inf \nu.$
 \end{corollary}
  \bpf Let
  $f = \phi ^{\frac{1}{n}}\nu^{\frac{1}{n}}(x)$
  and $h = \frac{1}{2}(1 + |\gra u|^2).$ We only need to check the
  conditions on $h.$ Choose $\Lambda = 1$ and $M = 1.$ It is easy
  to see that they satisfy the required conditions.
  \epf
 \vskip 1em

 3. \emph{Convex Hypersurface and Gauss Curvature Equation}
  \vskip 1em
   Given a closed smooth embedded $(n-1)$-dimensional submanifold $\Sigma$ in
   $\mathbb{R}^{n+1}$, we are asked whether there exists a hypersurface of
   constant Gauss-Kronecker curvature in $\mathbb{R}^{n+1}$ with $\Sigma$ as
   its boundary. Locally this problem is reduced to some
   Monge-Ampere type equation. If a hypersurface is locally
   strictly convex, we can express it as a graph $(x, u(x))$ for
   $x \in \Omega \subset \mathbb{R}^n$ which satisfies
  \beq \label{e:gauss1}
  \det^{\qquad \frac{1}{n}}(\hess u)= \kappa^{\frac{1}{n}} (x)
  ( 1 +|\gra u|^2)^{\frac{n+2}{2n}}
  \eeq
  where $\kappa$ is the Gauss curvature which is positive.
  In Caffarelli-Nirenberg-Spruck \cite{CNS-1}, this type of Monge-Ampere
  equation is studied in a strictly convex domain $\Omega.$
  On the other hand,  Guan-Spruck \cite{GbS93} consider
   star-shaped regions, i.e., radial graphs over a domain
  $\Omega \subset S^n$. In this setting, the problem becomes
  finding a positive solution $v(x)$ of the following equation in
  $\Omega \subset S^n:$
 $$\det^{\qquad \frac{1}{n}}(g_c ^{-1} (\hess v + v g_c))
  = \kappa^{\frac{1}{n}} (x) ( \frac{v^2 + |\gra v|^2}{v^2})^{\frac{n+2}{2n}}$$
  where $g_c$ is the standard metric on $S^n.$
  Since $v(x)$ is positive, let $u(x) = \ln v(x)$. The equation
  becomes
  \beq \label{e:gauss2}
  \det^{\qquad \frac{1}{n}}(g_c ^{-1} (\hess u + du \otimes du + g_c))
  = \kappa^{\frac{1}{n}} (x) ( 1+ |\gra u|^2)^{\frac{n+2}{2n}}
  e^{-u}.
  \eeq

  Equation (\ref{e:gauss2}) is in the form of
  (\ref{e:csc}) now. It is proved in Guan-Spruck \cite{GbS93}
  that the Hessian bound of the solution $u(x)$ to
  (\ref{e:gauss2}) over $\overline{\Omega}$ is less than or equal to that on $\de
  \Omega.$ Here, as a corollary of Theorem~\ref{t:max}, we have the
  following.

  \begin{corollary}
   Let $u(x)$ be a $C^4$ solution to (\ref{e:gauss1}) (or
   (\ref{e:gauss2})) in a bounded domain $\Omega \subset \mathbb{R}^n$ (or
   $S^n$, respectively) with $\hess u \in \Gamma^+_n$ (or $\hess u +
   du \otimes du + g_c \in \Gamma^+_n$, respectively).
   Suppose that $\kappa$ is positive.  Then
     $$ \sup_{\overline{\Omega}} |\hess u| < C $$
    where $C$ depends on $n, \|\kappa\|_{C^2}, \|u\|_{C^1(\overline{\Omega})}, \sup_{\de
    \Omega} |\hess u|$ and $\inf \kappa.$
  \end{corollary}
   \bpf For equation (\ref{e:gauss1}), $K=a=0$, and for equation
   (\ref{e:gauss2}), $K = a =1.$ In both cases, $h = (1 + |p|^2)^{\frac{n+2}{2n}}.$
    We only need to check the convexity condition on $h.$
    If $n=2$, then $h_{p_i p_j} = 2 \gd_{ij}.$  When $n> 2,$
    \begin{eqnarray*}
    h_{p_i p_j} &=& \frac{n+2}{n} (1+ |p|^2)^{-\frac{n-2}{2n}}
     \left( -\frac{n-2}{2n} \frac{p_i p_j}{1+ |p|^2} +
     \gd_{ij}\right)\\
      &\geq& \frac{n+2}{n} (1+ |p|^2)^{-\frac{n-2}{2n}}
      \left(\frac{1+ \frac{2}{n}|p|^2}{1+ |p|^2} \gd_{ij}\right)
        \geq \frac{n+2}{n} (1+ |p|^2)^{-\frac{3n-2}{2n}} \gd_{ij}
    \end{eqnarray*}
    Hence, $h_{p_i p_j} > \epsilon \gd_{ij}$ with $\epsilon$
    depending on $\sup |\gra u|.$
   \epf


 \section{Proof of Theorem~\ref{t:local}} \label{s:local}
 \bpf  We always assume $g_{ij} = \gd_{ij}$ at the point we are
 evaluating. Let $W = \hess u + a(x) du \otimes du + b(x) |\gra u|^2 g + B(x).$
 We will show that $\Delta u$ is bounded. By the condition $\Gamma \subset \Gamma^+_1,$
 we have
 $$0 < tr_g W = \Delta u + ( a(x)+ nb(x)) |\gra u|^2 + tr_g B(x). $$
 Since $a(x)+ nb(x) < - \gd_2$, the Laplacian $\Delta u$ has
 lower bound and
 \begin{equation} \label{i:grad1}
 |\gra u|^2 < C (\Delta u + 1 ).
 \end{equation} Therefore, we may assume $\Delta u$ is positive.
  Let $H = \eta (\Delta u + a(x)|\gra u| ^2) = \eta L $
 where $0 \leq \eta \leq 1$ is a cutoff function such that $\eta = 1$ in $B_{\frac{r}{2}}$
 and $\eta = 0$ outside $B_r$, and also $|\gra \eta| < C \frac{\sqrt{\eta}}{r}$
 and $|\hess \eta| < \frac{C}{r^2}.$ Without loss of generality,
 we assume $r=1$ since for general $r$, the proof is similar.

 Now by the condition $\Gamma \subset \Gamma^+_1$ again, we get
 \beq \label{i:l}
  L > -n b |\gra u|^2 - tr_g B \geq \gd_1 n |\gra u|^2 - tr_g B > -
  C.
 \eeq Hence, $L$ is lower bounded and we only need to get the upper bound
 of $L$. Suppose $x_0$ is the maximal point of $H$.
 At $x_0$, we have
 \begin{equation} \label{e:star1}
 H_i = \eta_i L + \eta L_i = \eta_i (\Delta u + a(x)|\gra u| ^2)+
 \eta (u_{kki} + a_i |\gra u|^2 + 2a u_k u_{ki}) = 0,
 \end{equation}
and
 $$ H_{ij} = \eta_{ij} L + \eta_i L_j + \eta_j L_i + \eta L_{ij}
    = (\eta_{ij} -2 \frac{\eta_i \eta_j}{\eta}) L + \eta
        L_{ij}
 $$ is negative semi-definite where in the second equality we have used
 (\ref{e:star1}). Moreover,
 $$L_{ij} = u_{kkij} + a_{ij} |\gra u|^2 + 2 a_i u_k u_{kj} + 2 a_j u_k u_{ki}
 + 2a u_{ki} u_{kj} + 2a u_k u_{kij}.$$
Using the positivity of $F^{ij}$ and the condition on $\eta$, we
get
 \beq \label{i:est}
  0 \geq F^{ij} H_{ij} = F^{ij}((\eta_{ij} -2 \frac{\eta_i \eta_j}{\eta}) L + \eta
        L_{ij}) \geq -C \sum_i F^{ii} L + \eta F^{ij} L_{ij}.
 \eeq
 Now to compute $ F^{ij} L_{ij}$, we note that $F^{ij} (2 a_i
u_k u_{kj}) = F^{ij}(2 a_j u_k u_{ki})$ because $F^{ij}$ is
symmetric. Thus, we obtain
 $$F^{ij} L_{ij} = F^{ij}(u_{kkij} + a_{ij} |\gra u|^2 + 4 a_i u_k u_{kj}
 + 2a u_{ki} u_{kj} + 2a u_k u_{kij}).$$
Changing the order of the covariant differentiations and using
(\ref{i:grad1}) give
\begin{eqnarray*}
 F^{ij} L_{ij} &\geq& F^{ij} u_{ijkk} + F^{ij}(2a u_{ki} u_{kj}+ 2a u_k u_{ijk})-
   C\sum_i F^{ii} ( 1 + |\hess u|^{\frac{3}{2}})\\
    &=& I + II - C\sum_i F^{ii}( 1 + |\hess u|^{\frac{3}{2}}).
 \end{eqnarray*}
 To compute I, notice that
 \begin{eqnarray*}
 W_{ij,kk} &=& u_{ijkk} + \Delta a u_i u_j + 2 a_k u_{ik} u_j + 2 a_k u_{jk}
 u_i + a (u_{ikk}u_j + 2 u_{ik} u_{jk} + u_i u_{jkk})\\
 & &+ (\Delta b |\gra u|^2 +
  4 b_k u_{kl} u_l + 2b |\hess u|^2 + 2b u_l u_{lkk}) \gd_{ij} +
  B_{ij,kk}.
 \end{eqnarray*} Then
 \begin{eqnarray*}
  I = F^{ij} u_{ijkk} &=& F^{ij}( W_{ij,kk}- \Delta a u_i u_j - 4 a_k u_{ik} u_j
  - 2a ( u_{ikk}u_j + u_{ik} u_{jk})\\
  & &- (\Delta b |\gra u|^2 +
  4 b_k u_{kl} u_l + 2b |\hess u|^2 + 2b u_l u_{lkk}) \gd_{ij} -
  B_{ij,kk})\\
  &\geq& F^{ij}W_{ij,kk}+ F^{ij}( - 2a ( u_{ikk}u_j +  u_{ik}
  u_{jk})- 2b (|\hess u|^2 + u_l u_{lkk}) \gd_{ij})\\
  & & - C\sum_i F^{ii}( 1 + |\hess u|^{\frac{3}{2}}).
 \end{eqnarray*}
 Changing the order of covariant differentiations again yields
 \begin{eqnarray*}
  I &\geq& F^{ij}W_{ij,kk}+ F^{ij}( - 2a ( u_{kki}u_j +  u_{ik}
  u_{jk})- 2b (|\hess u|^2 + u_l u_{kkl}) \gd_{ij})\\
  & & - C\sum_i F^{ii}( 1 + |\hess u|^{\frac{3}{2}}).
 \end{eqnarray*}
 Now we replace the terms $u_{kki}$ and $u_{kkl}$ by (\ref{e:star1}) to get
  \begin{eqnarray*}
  I &\geq& F^{ij}W_{ij,kk}+ F^{ij}( 2a u_j( a_i |\gra u|^2 + 2a u_k u_{ki}+
  \frac{\eta_i}{\eta}L ) - 2a u_{ik} u_{jk} - 2b |\hess u|^2 \gd_{ij} \\
  & & + 2b u_l ( a_l |\gra u|^2 + 2a u_k u_{kl}+
  \frac{\eta_l}{\eta}L ) \gd_{ij}) - C\sum_i F^{ii}( 1 + |\hess u|^{\frac{3}{2}}).
  \end{eqnarray*}
 Using (\ref{i:grad1}) again and the condition on $\eta,$ we have
  \begin{eqnarray*}
  I  &\geq& F^{ij}W_{ij,kk}+ F^{ij}( 4a^2 u_k u_{ki} u_j - 2a u_{ik}
   u_{jk} - 2b |\hess u|^2 \gd_{ij}+ 4ab u_l u_k u_{kl} \gd_{ij})\\
  & & - C \eta^{-\frac{1}{2}} \sum_i F^{ii} |\gra u| L
  - C\sum_i F^{ii}( 1 + |\hess u|^{\frac{3}{2}}).
 \end{eqnarray*}

 For II, we use the formula
 $$W_{ij,k} = u_{ijk} + a_k u_i u_j + au_{ik}u_j + a u_{jk}u_i
 + b_k |\gra u|^2 \gd_{ij} + 2b u_l u_{lk} \gd_{ij}+ B_{ij,k}$$
 to obtain
 \begin{eqnarray*}
  II &=& F^{ij}(2a u_{ki} u_{kj}+ 2a u_k u_{ijk})= F^{ij}(2a u_{ki}
  u_{kj} + 2a u_k W_{ij,k}\\
  & &+ 2a u_k (- a_k u_i u_j - 2au_{ik}u_j - b_k |\gra u|^2 \gd_{ij}
  - 2b u_l u_{lk} \gd_{ij}- B_{ij,k}))\\
  &\geq&  2a u_k F^{ij}W_{ij,k} + F^{ij}(2a u_{ki} u_{kj}- 4a^2 u_k u_{ik}u_j
 - 4ab u_l u_{lk} \gd_{ij})\\
  & &- C\sum_i F^{ii}( 1 + |\hess u|^{\frac{3}{2}}).
 \end{eqnarray*}
 Combining I and II together, we find that
 \begin{eqnarray*}
   F^{ij} L_{ij} &\geq& I + II - C\sum_i F^{ii}( 1 + |\hess u|^{\frac{3}{2}})\\
    &\geq&  F^{ij}W_{ij,kk} + 2a u_k F^{ij}W_{ij,k}+ F^{ij}( 4a^2 u_k u_{ki} u_j
    - 2a u_{ik} u_{jk} - 2b |\hess u|^2 \gd_{ij}\\
    & &+ 4ab u_l u_k u_{kl} \gd_{ij})
    + F^{ij}(2a u_{ki} u_{kj}- 4a^2 u_k u_{ik}u_j - 4ab u_l u_{lk} \gd_{ij})\\
  & & - C \eta^{-\frac{1}{2}} \sum_i F^{ii} |\gra u| L
  - C\sum_i F^{ii}( 1 + |\hess u|^{\frac{3}{2}}).
 \end{eqnarray*}
 After the cancellations, finally we arrive at
 \begin{eqnarray*}
   F^{ij} L_{ij}&\geq& F^{ij}W_{ij,kk} + 2a u_k F^{ij}W_{ij,k} - 2b \sum_i F^{ii}|\hess
   u|^2\\
   & & - C \eta^{-\frac{1}{2}} \sum_i F^{ii} |\gra u| L
    - C\sum_i F^{ii}( 1 + |\hess u|^{\frac{3}{2}})
 \end{eqnarray*}
 Now returning to (\ref{i:est}) and applying $\eta$ on both sides
 produces
 \begin{eqnarray*}
  0&\geq& \eta F^{ij} H_{ij} \geq -C \eta \sum_i F^{ii} L + \eta^2 F^{ij}
  L_{ij}\\
   &\geq & \eta^2 F^{ij}W_{ij,kk} + 2a \eta^2 u_k F^{ij}W_{ij,k}
   - 2b \eta^2 \sum_i F^{ii}|\hess u|^2 - C \eta^{\frac{3}{2}} \sum_i F^{ii}|\gra u| L \\
   & &-C \eta \sum_i F^{ii} L  - C \eta^2 \sum_i F^{ii}( 1 + |\hess u|^{\frac{3}{2}})\\
   &\geq& \eta^2 F^{ij}W_{ij,kk} + 2a \eta^2 u_k F^{ij}W_{ij,k}
   - 2b \eta^2 \sum_i F^{ii}|\hess u|^2 \\
   & & - C \sum_i F^{ii}( 1 + \eta |\hess u| + (\eta |\hess u|)^{\frac{3}{2}}).
 \end{eqnarray*}
 By the concavity of $F$ and Lemma~\ref{l:sym} (a), we have $\; F^{ij}W_{ij,kk}
 \geq (F^{ij}W_{ij})_{kk}\\ = (f(x,u) h(x, \gra u))_{kk}.$ Hence,
 \begin{align}
   0&\geq \eta^2 (f(x,u) h(x,\gra u))_{kk} + 2a \eta^2 u_k (f(x,u)h(x,\gra u))_k
   - 2b \eta^2 \sum_i F^{ii}|\hess u|^2 \notag\\
   &   - C \sum_i F^{ii}( 1 + \eta |\hess u|+ (\eta |\hess u|)^{\frac{3}{2}}). \label{i:fh}
 \end{align}

 case(a): $h$ is a positive constant. By Lemma~\ref{l:sym} (b), $\sum_i F^{ii} \geq F(e) = 1,$
    hence
   \begin{eqnarray*}
   0 &\geq& \eta^2 \frac{\de^2 f(x,u(x))}{\de x_k^2} h  + 2a \eta^2 u_k \frac{\de f(x,u(x))}{\de x_k} h
   - 2b \eta^2 \sum_i F^{ii}|\hess u|^2\\
   & &- C \sum_i F^{ii}( 1 + \eta |\hess u|+ (\eta |\hess u|)^{\frac{3}{2}})\\
   &\geq& - 2b \eta^2 \sum_i F^{ii}|\hess u|^2
   - C \sum_i F^{ii}( 1 + \eta |\hess u|+ (\eta |\hess u|)^{\frac{3}{2}}).
   \end{eqnarray*}
 By the condition on $b$, finally we arrive at
  $$
   0 \geq  \sum_i F^{ii} ( 2\gd_1 (\eta |\hess u|)^2  - C (\eta |\hess
    u|)- C (\eta |\hess u|)^{\frac{3}{2}} - C).$$
  This gives $(\eta |\hess u|) (x_0) \leq C$ and
  hence $H(x) =  \eta (\Delta u + a|\gra u| ^2)= (\Delta u + a|\gra u|
  ^2)= L$ is bounded for $x \in B_{\frac{1}{2}}.$ Now by (\ref{i:l}),
  we see $\gd_1 n |\gra u|^2 \leq L + tr_g B \leq C,$ which implies $|\gra u|^2$
   is bounded. And then $\Delta u = L - a |\gra u|^2$ is bounded.
 \vskip 1em

 case(b): $h = h(\gra u)$ and $f = f(x).$ First we perform some
 computations:
  \begin{eqnarray*}
   (f(x) h(\gra u))_{kk} &=& f_{kk} h + 2 f_k  h_{p_i} u_{ik} + f
  h_{p_i p_j} u_{ik} u_{jk} + f h_{p_i} \, u_{ikk}\\
    &\geq& - C \Lambda (1 + |\hess u|^{\frac{3}{2}}) + f \Lambda |\hess
    u|^2 + f h_{p_i} \, u_{ikk}
  \end{eqnarray*}
  where we have used the conditions on $h.$ Now changing the order of
  differentiations of $u_{ikk}$ and
   using (\ref{e:star1}) to replace $u_{kki}$ give
   \begin{eqnarray*}
    (f(x) h(\gra u))_{kk} &\geq& -C \Lambda (1 + |\hess u|^{\frac{3}{2}})
    + f \Lambda |\hess u|^2 - f h_{p_i} ( a_i |\gra u|^2 + 2 a u_k u_{ki} + \frac{\eta_i}{\eta}
    L)\\
     &\geq& -C \Lambda (1 + |\hess u|^{\frac{3}{2}})
    + f \Lambda |\hess u|^2 - 2a f h_{p_i} u_k u_{ki} -
    \frac{C}{\sqrt{\eta}} f |\gra_p h| L.
   \end{eqnarray*}
    On the other hand, we have
   $$2 a u_k (f(x)h(\gra u))_k = 2 a f_k h u_k  + 2 a f h_{p_i}\, u_{ik} u_k
      \geq -C \Lambda (1 + |\hess u|^{\frac{3}{2}}) + 2 a f h_{p_i}\, u_{ik} u_k.$$
   Thus returning to (\ref{i:fh}), we get
   \begin{eqnarray*}
   0 &\geq& -C \eta^2\Lambda (1 + |\hess u|^{\frac{3}{2}})
   + f \eta^2 \Lambda |\hess u|^2 - 2a \eta^2 f h_{p_i} u_k u_{ki} -
    C \eta^{\frac{3}{2}} f |\gra_p h| L\\
    & & + 2 a \eta^2 f h_{p_i}\, u_{ik} u_k - 2b \eta^2 \sum_i F^{ii}|\hess u|^2
     - C \sum_i F^{ii}( 1 + \eta |\hess u|+ (\eta |\hess
     u|)^{\frac{3}{2}})\\
    &\geq& -C \eta^{\frac{3}{2}} \Lambda (1 + |\hess u|^{\frac{3}{2}})
     + f \eta^2 \Lambda |\hess u|^2 - 2b \eta^2 \sum_i F^{ii}|\hess
     u|^2\\
    & & - C \sum_i F^{ii}( 1 + \eta |\hess u|+ (\eta |\hess
     u|)^{\frac{3}{2}}).
   \end{eqnarray*}
  Applying the conditions $b$ and using Lemma~\ref{l:sym} (b) to obtain
    \begin{eqnarray*}
   0 &\geq&  \Lambda ( - C - C \eta^{\frac{3}{2}} |\hess u|^{\frac{3}{2}}
    + c_{inf}\eta^2 |\hess u|^2)\\
    & & + \sum_i F^{ii} (2\gd_1 \eta^2 |\hess u|^2
     - C  -C \eta |\hess u| -C \eta^{\frac{3}{2}} |\gra u||\hess
     u|)\\
    &\geq&  -C + \sum_i F^{ii} (\gd_1 \eta^2 |\hess u|^2 - C )
    \geq  \sum_i F^{ii} (\gd_1 \eta^2 |\hess u|^2 - C) .
   \end{eqnarray*}

    This gives $(\eta |\hess u|) (x_0) \leq C$ and then
  $H \leq C.$ Therefore, $\Delta u$ and $|\gra u|^2$ are all bounded.
   \vskip 1em
  case(c): $|\gra u|$ is bounded and thus $h$ is bounded. This gives
   $$(f(x,u) h(x,\gra u))_{kk} + 2a u_k (f(x,u)h(x,\gra u))_k \geq
   -C (1 + |\hess u|^2 ) + f h_{pi} u_{ikk}.$$
   We change the order of differentiations of third derivative terms and
   use (\ref{e:star1}) to replace $u_{kki}:$
   $$(f(x,u) h(x,\gra u))_{kk} + 2a u_k (f(x,u)h(x,\gra u))_k \geq
     -C (1 + |\hess u|^2 ) - \frac{C}{\sqrt{\eta}}(1 + |\hess u|). $$
   Hence, (\ref{i:fh}) becomes
   \begin{eqnarray*}
    0&\geq&  -C \eta^2(1 + |\hess u|^2 ) - C \eta^{\frac{3}{2}}(1 + |\hess
    u|))
    +\sum_i F^{ii}(- 2b \eta^2|\hess u|^2 -C - C \eta |\hess u|)\\
    &\geq& -C - C \eta^2 |\hess u|^2 +\sum_i F^{ii}(- 2b \eta^2|\hess u|^2 -C - C \eta |\hess
    u|).
   \end{eqnarray*}
   By (A) and condition on $b$, we see that
    \begin{eqnarray*}
    0&\geq& -C - C \eta^2 |\hess u|^2 + \sum_i F^{ii} (2\gd_1 \eta^2|\hess u|^2
    - C  - C \eta |\hess u|)\\
     &\geq & -C - C \eta^2 |\hess u|^2 + \mu_0(\frac{\gs_1}{F})^{\mu_1}
     (2\gd_1 \eta^2|\hess u|^2- C).
    \end{eqnarray*}
   Apply $(\eta F)^{\mu_1}$ on both sides and
   note that $\gs_1 = \Delta u + ( a(x)+ nb(x)) |\gra u|^2 + tr\, B(x) \geq
   \Delta u -C ,$ so we have
     \begin{eqnarray*}
     0 &\geq & -C - C \eta^2 |\hess u|^2 + \mu_0 \gs_1^{\mu_1}
     (2\gd_1 \eta^{2+ \mu_1}|\hess u|^2- C\eta^{\mu_1})\\
       &\geq& -C - C \eta^2 |\hess u|^2 + 2\gd_1 \mu_0 \eta^{2+ \mu_1}(\Delta u)^{\mu_1}
     |\hess u|^2- C\eta^{\mu_1}(\Delta u)^{\mu_1}.
     \end{eqnarray*}
      This gives $(\eta \Delta u) (x_0) \leq C$, and consequently $\Delta u$ is bounded.
 \vskip 1em
  Once $\Delta u$ is bounded, to get the Hessian bounds for cases (a) and (b), we simply
  consider the maximum of the tensor $\eta (\hess u + a du \otimes du)$ over the set
  $(x, \xi) \in (B_1, S^n)$. As for case (c), we use the basic
  fact that if $\Gamma^+_2 \subset \Gamma,$ then $- \frac{n-2}{n} \gs_1 \leq \gl_i \leq \gs_1$
  for $\gl \in \Gamma.$
 \epf
\vskip 1em

\section{Proof of Theorem~\ref{t:max}}\label{s:max}
\bpf We assume $g_{ij} = \gd_{ij}$ at the point we are
 evaluating. Now we start with some computations on curvatures.
 It is known that the Riemannian curvature has the
 decomposition
 $$R_{ijkl} = \mathcal{W}_{ijkl} + (A_{ik}g_{jl} + A_{jl}g_{ik}- A_{il}g_{jk} - A_{jk}g_{il}) $$
 where $\mathcal{W}$ is the Weyl tensor and $A$ is the Schouten
 tensor. If $g$ is of constant sectional curvature $K$, then $\mathcal{W}$ is zero,
 $Ric = (n-1) K g$ and $R = n(n-1)K. $ Hence we have
 $$ R_{ijkl} = K (g_{ik}g_{jl}- g_{il}g_{jk}).$$

Let $W = \hess u + a du \otimes du+ Kg.$  By $\Gamma \subset
\Gamma^+_1,$ we get $\Delta u + a |\gra u|^2 + n K > 0.$ Thus
$\Delta u$ is lower bounded and we only need to get the upper
 bound. Let $H = \Delta u + a|\gra u| ^2 .$ We may assume $H$ is
 large and suppose $x_0$ is the maximal point of $H$.
 At $x_0$, we have
 \begin{equation} \label{e:star2}
 H_i = u_{kki} + 2a u_k u_{ki} = 0,
 \end{equation}
and
$$ H_{ij} = u_{kkij} + 2a u_{ki} u_{kj} + 2a u_k u_{kij}$$ is negative
semi-definite.
Using the positivity of $F^{ij}$ ,we get
 $$ 0 \geq F^{ij} H_{ij} = F^{ij}u_{kkij}+ F^{ij}(2a u_{ki} u_{kj} +
2a u_k u_{kij})= I + II.$$
 Before computing I and II, we examine carefully the formulae at the end of
 Section~\ref{s:sigma}:
 \begin{align}
 u_{kij} &= u_{ijk} + R_{mikj} u_m = u_{ijk} + K(g_{ij} u_k - g_{ik}
 u_j),\notag \\
 u_{kkij} &= u_{ijkk} + 2 R_{mikj} u_{mk} - R_{mj} u_{mi} - R_{mi} u_{mj}=
 u_{ijkk} + 2K \Delta u g_{ij} -2Kn u_{ij}.\notag
 \end{align}
 Thus I becomes
  $$I = F^{ij} u_{kkij} = F^{ij} (u_{ijkk} + 2 K \Delta u g_{ij} - 2n K u_{ij}). $$
 Now use the formula
 $$ W_{ij,kk} = u_{ijkk} + a (u_{ikk}u_j + 2 u_{ik} u_{jk} + u_i
 u_{jkk})$$
 to get
 $$ I = F^{ij}( W_{ij,kk}- 2a ( u_{ikk}u_j + u_{ik}
 u_{jk}) + 2K \Delta u g_{ij} -2Kn u_{ij}),$$
 where we have used $F^{ij} a u_{ikk}u_j= F^{ij} a u_{jkk}u_i$ because $F^{ij}$
  is symmetric. Changing the order of the differentiations of $u_{ikk}$ and
  replacing it by (\ref{e:star2}) gives
  \begin{eqnarray*}
 I &=& F^{ij}( W_{ij,kk}- 2a(u_{kki} + (n-1)K u_i)u_j -2a u_{ik}
 u_{jk} + 2K \Delta u g_{ij} -2Kn u_{ij})\\
   &=& F^{ij}( W_{ij,kk}+ 4a^2 u_{k}u_{ki}u_j - 2a(n-1)K u_iu_j -2a u_{ik}
 u_{jk} + 2K \Delta u g_{ij} -2Kn u_{ij}).
  \end{eqnarray*}
 For II, we first change the order of differentiations of $u_{kij}$
 and then replace $u_{ijk}$ by $W_{ij,k} - a u_i u_{jk} - a u_j u_{ik}$
 to get
 \begin{eqnarray*}
  II &=& F^{ij}(2a u_{ki} u_{kj}+ 2a u_k u_{kij})
    = F^{ij}(2a u_{ki} u_{kj}+ 2a u_k u_{ijk} + 2aK u_k (g_{ij} u_k - g_{ik}
    u_j))\\
    &=& F^{ij}(2a u_{ki} u_{kj} + 2a u_k W_{ij,k}- 4a^2 u_k u_{ik}u_j+
    2aK |\gra u|^2 g_{ij}- 2aK u_i u_j).
 \end{eqnarray*}

 We combine I and II, and note the cancellation. We obtain
 \begin{eqnarray*}
  I + II &=& F^{ij}W_{ij,kk} + 2a u_k F^{ij}W_{ij,k} + F^{ij} (- 2nK au_iu_j
     -2Kn u_{ij}+ 2aK |\gra u|^2 g_{ij}\\
     & & + 2K \Delta u g_{ij}).
 \end{eqnarray*}
  Replacing $u_{ij}$ by $W_{ij}- a u_i u_j - K g_{ij}$ and using the concavity of $F$,
   we get
  \begin{eqnarray*}
   0 &\geq& I + II = F^{ij}W_{ij,kk} + 2a u_k F^{ij}W_{ij,k} + F^{ij}(
     -2Kn W_{ij} + 2K^2n g_{ij} + 2K H g_{ij})\\
    &\geq& (f(x,u) h(\gra u))_{kk} + 2a u_k (f(x,u)h(\gra u))_k
     -2Kn f(x,u)h(\gra u) \\
    & & +  2K\sum_i F^{ii} (K n+ H ).
   \end{eqnarray*}
 Since we have $C^1$ bounds and nonnegative $K$, we obtain
  \begin{eqnarray*}
   0&\geq& (f(x,u) h(\gra u))_{kk} + 2a u_k (f(x,u)h(\gra u))_k
     -2Kn f(x,u)h(\gra u) \\
    &\geq& -C - C |\hess u| + f h_{p_i p_j} u_{ik} u_{jk} + f h_{p_i}
    u_{ikk}.
 \end{eqnarray*}
 Changing the order of the differentiations of $u_{ikk}$ and replacing it by (\ref{e:star2})
 again, produce
  \begin{eqnarray*}
   0 &\geq& -C - C |\hess u| + f h_{p_i p_j} u_{ik} u_{jk} +
    f h_{p_i} ( -2a u_k u_{ki} + (n-1)K u_i)\\
     &\geq& -C - C |\hess u| + f h_{p_i p_j} u_{ik} u_{jk}.
  \end{eqnarray*}
  Then by the convexity of $h$, we arrive at
   $$0 \geq -C - C |\hess u| + f \epsilon |\hess u| \geq -C - C |\hess u| +
   \epsilon \, c_{inf}|\hess u|^2.$$
   This gives $|\hess u(x_0)| < C$ and hence $H < C.$
   Finally, to get the Hessian bounds, we
  consider the tensor $\hess u + a du \otimes du$ over the set
  $(x, \xi) \in (\Omega, S^n)$.
 \epf


 Princeton University, Princeton, NJ \par
 Email address: \textsf{szuchen@math.princeton.edu}

\end{document}